\documentclass{amsart}
\usepackage{amssymb,latexsym,graphics,amsfonts}

\swapnumbers

\theoremstyle{definition}

\theoremstyle{definition}

\numberwithin{equation}{section}
 \numberwithin{equation}{subsection}

\begin{document}

\title[module topological center]{module amenability of the second dual and module topological center of semigroup algebras}

\author[M. Amini]{Massoud Amini}
\address{Department of Mathematics,
 Tarbiat Modares University , Tehran 14115-175, Iran}
\email{mamini@modares.ac.ir}

\author[A. Bodaghi]{Abasalt Bodaghi}
\address{Islamic AZAD University, Garmsar Branch, Garmsar,
 Iran.}
\email{abasalt\_bodaghi@yahoo.com}

\author[D. Ebrahimi Bagha]{Davood Ebrahimi Bagha}
\address{Faculty of Mathematical Sciences,
 Shahid Beheshti University, Evin, Tehran, Iran}
\email{d-ebrahimi@cc.sbu.ac.ir}

\address{current address of the first author: Institut Penyelidikan Matematik, universiti Putra Malaysia, 43400 UPM Serdang,
Selangor Darul Ehsan, Malaysia}

\email{massoud@putra.upm.edu.my}

\subjclass[2000]{46H25}

\keywords{Banach modules, module derivation, module amenability,
module topological center, inverse semigroup}

\dedicatory{}

\smallskip

\begin{abstract}
Let $S$ be an inverse semigroup with an upward directed set of
idempotents $E$. In this paper we define the module topological
center of second dual of a Banach algebra which is a Banach
module over another Banach algebra with compatible actions, and
find it for $ \ell ^{1}(S)^{**}$ (as an $\ell^{1}(E)$-module). We
also prove that $ \ell ^{1}(S)^{**}$ is $\ell^{1}(E)$-module
amenable if and only if an appropriate group homomorphic image of
$S$ is finite.
\end{abstract}

\maketitle

\section{introduction}

The first author in \cite{am1} introduced the concept of module
amenability and showed that for an inverse semigroup $S$, the
semigroup algebra $\ell^1(S)$ is module amenable as a Banach
module on $\ell^1(E)$, where $E$ is the set of idempotents of
$S$, if and only if $S$ is amenable (see also [2]). The first and
third authors showed in \cite{ame} that $\ell^1(S)$ is weak
module amenable, for each commutative inverse semigroup $S$. In
\cite{ra1, ra2} the concept of Arens module regularity is
introduced and it is shown that when $S$ is an inverse semigroup
with totally ordered subsemigroup $E$ of idempotents, then
$A=\ell^1(S)$ is module Arens regular if and only if an
appropriate group homomorphic image $S/\approx$ of $S$ is finite.
When $S$ is a discrete group, we have $S/\approx=S$.

In part two of this paper, we define the module topological center
of second dual $\mathcal A^{**}$ of a Banach algebra $\mathcal A$
which is a Banach $\mathfrak A$-module with compatible actions on
another Banach algebra $\mathfrak A$. We show that if an inverse
semigroup $S$ has an upward directed set of idempotents $E$, for
the semigroup algebra $\ell^1(S)$ as an $\ell^1(E)$-module, the
module topological center of $ \ell ^{1}(S)^{**}$ is $ \ell
^{1}(S/\approx)$. This could be considered as the module version
(for inverse semigroups) of a result of Lau and Losert \cite{llo}
which asserts that for any locally compact group $G$, the
topological center $L^1(G)^{**}$ is the same as $L^1(G)$, a fact
which is also proved (using a different proof) by Lau and Ulger
in \cite{lul}. The existing semigroup versions of this result
usually assume cancellation type properties. For instance, in
[13], Lau showed that for any discrete weakly cancellative
semigroup $S$, the topological center of $\ell^1(S)^{**}$ is
$\ell^1(S)$. For the locally compact case, it is shown by Bami in
\cite{las} that for a locally compact, Housdorff, cancellative,
foundation topological *-semigroup with identity, such that
$C^{-1}D$ is compact for any compact subsets $C$ and $D$ of $S$,
we have $Z_t(M_a(S)^{**})=M_a(S)$, where $M_a(S)$ denotes the
space of all measures $ \mu \in M(S)$ (the space of all bounded
complex Radon measure on $S$) for which the mappings $x \mapsto
|\mu|*\delta_x $ and $x \mapsto \delta_x*|\mu|$ from $S$ into
$M(S)$ are weakly continuous. There are similar results by Filali
and Salmi \cite{fis} for the Beurling algebra of a weakly
cancellative, right cancellative, discrete semigroup with a
diagonally bounded weight.

In part three, under some mild conditions, we show that module
amenability of the second dual Banach algebra implies the module
amenability of the algebra. we show that $ \ell ^{1}(S)^{**}$ is
$\ell^{1}(E)$-module amenable if and only if $S/\approx$ is
finite. For the bicyclic semigroup $\mathcal C$, we show that
$\mathcal C/\approx\simeq\mathbb Z$ and conclude that $ \ell
^{1}(\mathcal C)^{**}$ is not module amenable. Since $\mathcal C$
is an amenable semigroup, we already know that $\ell
^{1}(\mathcal C)$ is module amenable. The fact that amenability
of $A^{**}$ implies amenability of $A$ is first proved by
Gourdeau in \cite{gou}. Different proofs are provided by
Ghahramani, Loy and Willis in \cite{glw}. Also it was first proved
in \cite{glw} that for a locally compact group $G$, the
amenability of $L^1(G)^{**}$ implies that $G$ is finite. We are
not aware of any similar result for semigroups.

\section{Module Topological Center }

Throughout this paper, ${\mathcal A}$ and ${\mathfrak A}$ are
Banach algebras such that ${\mathcal A}$ is a Banach ${\mathfrak
A}$-bimodule with compatible actions, that is

$$
\alpha\cdot(ab)=(\alpha\cdot a)b,
\,\,(ab)\cdot\alpha=a(b\cdot\alpha) \hspace{0.3cm}(a,b \in
{\mathcal A},\alpha\in {\mathfrak A}).
$$
Let ${\mathcal X}$ be a Banach ${\mathcal A}$-bimodule and a
 Banach ${\mathfrak A}$-bimodule with compatible actions, that is
$$
\alpha\cdot(a\cdot x)=(\alpha\cdot a)\cdot x,
\,\,a\cdot(\alpha\cdot x)=(a\cdot\alpha)\cdot x, \,\,(\alpha\cdot
x)\cdot a=\alpha\cdot(x\cdot a) \hspace{0.3cm}(a \in{\mathcal
A},\alpha\in {\mathfrak A},x\in{\mathcal X} )
$$
and the same for the right or two-sided actions. Then we say that
${\mathcal X}$ is a Banach ${\mathcal A}$-${\mathfrak A}$-module.
If moreover
$$\alpha\cdot x=x\cdot\alpha \hspace{0.3cm}( \alpha\in {\mathfrak
A},x\in{\mathcal X} )$$ then $\mathcal X $ is called a {\it
commutative} ${\mathcal A}$-${\mathfrak A}$-module. If $\mathcal X
$ is a (commutative) Banach ${\mathcal A}$-${\mathfrak
A}$-module, then so is $\mathcal X^*$, where the actions of
$\mathcal A$ and ${\mathfrak A}$ on $\mathcal X^*$ are defined by
$$\langle\alpha\cdot f,x\rangle{}=\langle{}f,x\cdot\alpha\rangle{},\,\,\langle{}
a\cdot f,x\rangle{}=\langle{}f,x\cdot a\rangle{}\hspace{0.3cm} (a
\in{\mathcal A},\alpha\in {\mathfrak A},x\in{\mathcal X},f \in
\mathcal X^* )$$
 and the same for the right actions. Let ${\mathcal Y}$ be another ${\mathcal A}$-${\mathfrak
A}$-module, then a ${\mathcal A}$-${\mathfrak A}$-module morphism
from ${\mathcal X}$ to ${\mathcal Y}$ is a norm-continuous map
$\varphi :{\mathcal X}\longrightarrow {\mathcal Y}$ with $
\varphi(x \pm y)= \varphi (x)\pm
 \varphi (y)$ and
 $$ \varphi (\alpha \cdot x)= \alpha \cdot\varphi(x), \,\,\varphi(x\cdot \alpha)= \varphi (x)\cdot\alpha
,\,\,\varphi (a\cdot x)=a\cdot \varphi(x), \varphi(x\cdot a)=
\varphi(x)\cdot a,$$ for $x,y \in {\mathcal X}, a \in {\mathcal
A},$ and $\alpha \in {\mathfrak A}$.

 Note that when
${\mathcal A}$ acts on itself by algebra multiplication, it is
not in general a Banach ${\mathcal A}$-${\mathfrak A}$-module, as
we have not assumed the compatibility condition
$$a\cdot(\alpha\cdot b)=(a\cdot\alpha)\cdot b\quad (\alpha\in {\mathfrak A}, a,b \in{\mathcal
A}).$$ If $\mathcal A$ is a commutative $\mathfrak A$-module and
acts on itself by multiplication from both sides, then it is also
a Banach ${\mathcal A}$-${\mathfrak A}$-module.

If ${\mathcal A}$ is a Banach $\mathfrak A$-module with
compatible actions, then so are the dual space $\mathcal A^*$ and
the second dual space $\mathcal A^{**}$. If moreover $\mathcal A$
is a commutative $\mathfrak A$-module, then $\mathcal A^*$ and the
$\mathcal A^{**}$ are commutative ${\mathcal A}$-${\mathfrak
A}$-modules. Also the canonical embedding $\hat{}: \mathcal A\to
\mathcal A^{**}; a\mapsto \hat a$ is an $\mathfrak A$-module
morphism.

Consider the projective tensor product $\mathcal A \widehat
\bigotimes \mathcal A$. It is well known that $\mathcal A
\widehat \bigotimes \mathcal A$ is a Banach algebra with respect
to the canonical multiplication map defined by
$$(a\otimes b)(c\otimes d)=(ac\otimes bd)$$
and extended by bi-linearity and continuity [4]. Then $\mathcal A
\widehat \bigotimes \mathcal A$ is a Banach ${\mathcal
A}$-${\mathfrak A}$-module with canonical actions. Let $I$ be the
closed ideal of the projective tensor product $\mathcal A
\widehat \bigotimes \mathcal A$ generated by elements of the form
$\alpha \cdot a \otimes b-a \otimes b\cdot\alpha$ for $ \alpha\in
{\mathfrak A},a,b\in{\mathcal A}$. Consider the map $\omega \in
{\mathcal L}(\mathcal A \widehat \bigotimes \mathcal A, \mathcal
A )$ defined by $\omega (a \otimes b)=ab$ and extended by
linearity and continuity. Let $J$ be the closed ideal of
${\mathcal A}$ generated by $\omega(I)$. Then the module
projective tensor product ${\mathcal A}\widehat \bigotimes
_{\mathfrak A} {\mathcal A}\cong(\mathcal A \widehat \bigotimes
\mathcal A)/{I}$ and the quotient Banach algebra $\mathcal A/J$
are Banach ${\mathfrak A}$-modules with compatible actions. Also
the map $\widetilde{\omega} \in {\mathcal L}({\mathcal A}\widehat
\bigotimes _{\mathfrak A} {\mathcal A},\mathcal A/J)$ defined by
$ \widetilde{\omega} (a \otimes b +I)=ab+J$ extends to an
${\mathfrak A}$-module morphism.

Let $ \square$ and $ \lozenge$ be the first and second Arens
products on the second dual space $\mathcal A^{**}$, then
$\mathcal A^{**}$ is a Banach algebra with respect to both of
these products. When these two products coincide on $\mathcal
A^{**}$, we say that $\mathcal A$ is Arens regular. Let
${\mathcal Z}_t(\mathcal A^{**})$ denote the topological center
of $\mathcal A^{**}$, that is
$${\mathcal Z}_t(\mathcal A^{**})=\{G \in \mathcal A^{**} :F
\mapsto G \square F \hspace{0.2cm} \text{is} \hspace{0.2cm}
 \sigma (\mathcal A^{**},\mathcal A^*)\text{-continuous} \}.$$
We define the {\it module topological center} of $\mathcal A^{**}$
(as an ${\mathfrak A}$-module) by
$$
{\mathcal Z}_{\mathfrak A}(\mathcal A^{**})=\{G \in \mathcal
A^{**} :F \longrightarrow G \square F \hspace{0.3cm} \text{is}
\hspace{0.2cm}
 \sigma (\mathcal A^{**},J^{\perp})\text{-continuous} \}.$$
Clearly ${\mathcal Z}_{\mathfrak A}(\mathcal A^{**})$ is a $
\sigma (\mathcal A^{**},J^{\perp})$-closed subalgebra of $
(\mathcal A^{**},\square)$ containing ${\mathcal A}$. Indeed
$\mathcal A\subseteq {\mathcal Z}_t(\mathcal A^{**}) \subseteq
{\mathcal Z}_{\mathfrak A}(\mathcal A^{**})$.

The Arens regularity of $\mathcal A$ is equivalent to the weak
compactness of linear maps $\mathcal R_\lambda : \mathcal A\to
\mathcal A^{*}; a \mapsto a\cdot\lambda$, for each $\lambda\in
\mathcal A^*$ (see \cite[Theorem 2.6.17]{dal} for more details).

\vspace{.3cm}\paragraph{\large\bf Definition 2.1.} \cite{ra2} {\it
$\mathcal A$ is called module Arens regular (as an ${\mathfrak
A}$-module) if ${\mathfrak A}$-module homomorphisms $R_\lambda$
are weakly compact for any $\lambda \in J^{\perp}$ satisfying
$\lambda(\alpha\cdot ab) = \lambda(ab\cdot\alpha),
\lambda(c(\alpha\cdot ab)) = \lambda(c(ab\cdot\alpha)),
\lambda((\alpha\cdot ab)d) = \lambda((ab\cdot\alpha)d)$, and
$\lambda(c(\alpha\cdot ab)d) = \lambda(c(ab\cdot\alpha)d)$, for
$\alpha\in {\mathfrak A}$ and $a, b, c, d \in{\mathcal A}$.}

\vspace{.3cm}\paragraph{\large\bf Proposition 2.2.} {\it The
following statements are equivalent:

(i) $\mathcal A$ is module Arens regular.

(ii) $ F \square G-F \lozenge G \in J^{\perp \perp}
,\hspace{0.3cm}F,G \in A^{**} $.

(iii) $ {\mathcal Z}_{\mathfrak A}(\mathcal A^{**})=\mathcal
A^{**}.$}

\vspace{.3cm}\paragraph{\large\bf Proof.} The equivalence of (i)
and (ii) is proved in \cite[Theorem 2.2]{ra1}.

(ii) $\Rightarrow$ (iii): Let $F,G \in \mathcal A^{**}$ and $G_j
\stackrel{J^{\perp}}{\longrightarrow}G$ (where the superscript
$J^\perp$ shows the convergence in the weak topology $\sigma
(\mathcal A^{**},J^\perp)$ generated the family $J^\perp$ of
$w^*$-continuous functionals on $\mathcal A^{**}$). Then $F
\lozenge G_j \stackrel{J^{\perp}}{\longrightarrow}F \lozenge G$
and $(F \square G_j-F \lozenge G_j)(f)=0,$ for any $j$ and $ f
\in J^{\perp}$, hence $ F \square G_j
\stackrel{J^{\perp}}{\longrightarrow}F \square G, $ therefore $F
\in {\mathcal Z}_{\mathfrak A}(\mathcal A^{**})$.

(iii) $\Rightarrow$ (ii): For $F,G \in \mathcal A^{**}$ there is a
bounded net $(b_k)\subset \mathcal A $ with $ \widehat b_k
\stackrel{J^{\perp}}{\longrightarrow}G$. Then $ F \square b_k
\stackrel{J^{\perp}}{\longrightarrow}F \square G, $ we have $ F
\square b_k=F \lozenge b_k$, hence $F \lozenge b_k
\stackrel{J^{\perp}}{\longrightarrow}F \square G, $ and since we
also have $F \lozenge b_k \stackrel{J^{\perp}}{\longrightarrow}F
\lozenge G, $ we get $F \square G=F \lozenge G$ on $ J^{\perp}
$.$\hfill\blacksquare $

\vspace{.3cm} By the proof of the above proposition, we have
$$ {\mathcal Z}_{\mathfrak A}(\mathcal A^{**})=\{G \in A^{**} : G \square
F-G \lozenge F \in J^{\perp \perp} \hspace{0.1cm} (F \in \mathcal
A^{**}) \}.$$

\vspace{.3cm}\paragraph{\large\bf Definition 2.3.} {\it A discrete
semigroup $S$ is called an inverse semigroup if for each $ s \in
S $ there is a unique element $s^* \in S$ such that $ss^*s=s$ and
$s^*ss^*=s^*$. an element $e \in S$ is called an idempotent if
$e=e^*=e^2$. The set of idempotents of $S$ is denoted by $E$.}

Throughout this section, $S$ is an inverse semigroup with set of
idempotents $E$, where the order of $E$ is defined by
$$e\leq d \Longleftrightarrow ed=e \hspace{0.3cm}(e,d \in
E).$$ Then $E$ is a commutative subsemigroup of $S$, and $ \ell
^{1}(E)$ could be regard as a subalgebra of $ \ell ^{1}(S)$, and
thereby $ \ell ^{1}(S)$ is a Banach algebra and a Banach $ \ell
^{1}(E)$-module with compatible actions. Here we let $ \ell
^{1}(E)$ act on $ \ell ^{1}(S)$ by multiplication from right and
trivially from left, that is
$$\delta_e\cdot\delta_s = \delta_s, \,\,\delta_s\cdot\delta_e = \delta_{se} =
\delta_s * \delta_e \hspace{0.3cm}(s \in S,  e \in E).$$

In this case, $J$ is the closed linear span of
$$\{\delta_{set}-\delta_{st} \hspace{0.2cm} s,t \in S,  e \in E
\}.$$ We consider an equivalence relation on $S$ as follows
$$s\approx t \Longleftrightarrow \delta_s-\delta_t \in J \hspace{0.2cm} (s,t \in
S).$$

It is shown in \cite[Theorem 3.2]{ra1} that if $E$ is totally
ordered, then the quotient $S/\approx$ is a discrete group. This
is a rather strong condition. Let us observe that the quotient
$S/\approx$ is a group under the weaker condition that $E$ is
upward directed. Recall that $E$ is called {\it upward directed }
if for every $e,f \in E$ there exist $g \in E$ such that $eg=e$
and $fg=f$. This is precisely the assertion that $S$ satisfies the
$D_1$ condition of Duncan and Namioka \cite{dna}. The bicyclic
semigroup is an example of an inverse semigroup with a totally
ordered set of idempotents. On the other hand, the free unital
inverse semigroup on two generators has an upward directed set of
idempotents which is not totally ordered, and finally the set of
idempotents of the free inverse semigroup on two generators (say
$a$ and $b$) is not even upward directed (as there is no
idempotent majorizing  both $aa^*$ and $bb^*$). If $E$ is
directed upward and $e,f \in E$, then
$\delta_g-\delta_f=\delta_{gg}-\delta_{gfg} \in J$, and so
$g\approx f$. Similarly $g\approx e$, hence $e\approx f$. Now the
argument of \cite[Theorem 3.2]{ra1} could be adapted to show that
in this case $S/\approx$ is again a discrete group.

\vspace{.3cm}\paragraph{\large\bf Theorem 2.4.}  {\it If $E$ is
upward directed then $Z_{\ell ^{1}(E)}(\ell ^{1}(S)^{**}) = {\ell
^{1}}(S/\approx)$.}

\vspace{.3cm}\paragraph{\large\bf Proof.} As in \cite[Theorem
3.3]{ra1}, we may observe that $\ell ^{1}(S)/J\cong {\ell
^{1}}(S/\approx)$. Also by the proof of \cite[Theorem 2.4]{ra1},
we have $G \square F-G \lozenge F \in J^{\perp \perp}$ if and
only if the images of $G \square F$ and $G \lozenge F$ in $(\ell
^{1}(S)^{**}/J^{\perp \perp},\square)$ and $(\ell
^{1}(S)^{**}/J^{\perp \perp},\lozenge)$ are equal. Therefore
\begin{align*}
Z_{\ell ^{1}(E)}(\ell ^{1}(S)^{**}) &=\{G \in \mathcal A^{**} : G
\square F-G \lozenge F \in J^{\perp \perp}, \hspace{0.3cm}
\forall F \in\ell ^{1}(S)^{**} \}  \\
&= \{G \in \mathcal A^{**} : G \square F=G \lozenge F
\hspace{0.2cm} \text{on} \hspace{0.2cm} \ell
^{1}(S)^{**}/J^{\perp \perp} \} \\
&=Z_t((\ell ^{1}(S)/J)^{**}) = Z_t({\ell ^{1}}(S/\approx)^{**}) =
{\ell ^{1}}(S/\approx).
\end{align*}
The last equality follows from \cite[Theorem
1]{llo}.$\hfill\blacksquare $

\section{Module Amenability  }

Let ${\mathcal A}$ and ${\mathfrak A}$ be as in the above section
and ${\mathcal X}$ be a Banach ${\mathcal A}$-${\mathfrak
A}$-module. Let $I$ and $J$ be the corresponding closed ideals of
$\mathcal A \widehat \bigotimes \mathcal A$ and $\mathcal A$,
respectively. A bounded map $D: \mathcal A \longrightarrow
\mathcal X $ is called a {\it module derivation} if
$$D(a\pm b)=D(a)\pm D(b),\hspace{0.2cm}D(ab)=D(a)\cdot b+a\cdot D(b)\hspace{0.3cm}
(a,b \in \mathcal A),$$and
$$D(\alpha\cdot a)=\alpha\cdot D(a),\hspace{0.3cm}D(a\cdot\alpha)=D(a)\cdot\alpha
\hspace{0.3cm}(a \in{\mathcal A},\alpha\in {\mathfrak A}).$$
Although $D$ is not necessary linear, but still its boundedness
implies its norm continuity (since it preserves subtraction).
When $\mathcal X $ is commutative, each $x \in \mathcal X $
defines a module derivation
$$D_x(a)=a\cdot x-x\cdot a \hspace{0.3cm} (a \in{\mathcal A}).$$
These are called {\it inner} module derivations. The Banach
algebra ${\mathcal A}$ is called {\it module amenable} (as an
${\mathfrak A}$-module) if for any commutative Banach ${\mathcal
A}$-${\mathfrak A}$-module $\mathcal X $, each module derivation
$D: \mathcal A \longrightarrow \mathcal X^*$ is inner \cite{am1}.

\vspace{.3cm}\paragraph{\large\bf Lemma 3.1.} {\it Let $\mathcal
A$ be a Banach algebra and Banach ${\mathfrak A}$-module with
compatible actions, and $J_0$ be a closed ideal of $\mathcal A$
such that $J \subseteq J_0$. If $\mathcal A/J_0$ has a left or
right identity $e+J_0$, then  for each $\alpha \in {\mathfrak A}$
and $a\in\mathcal A$ we have $a\cdot\alpha-\alpha\cdot a \in
J_0$, i.e., $\mathcal A/J_0$ is commutative Banach ${\mathfrak
A}$-module.}

\vspace{.3cm}\paragraph{\large\bf Proof.} We prove the result for
the left identity. For each $\alpha \in {\mathfrak A}$ and
$a\in\mathcal A$, we have
$(e+J_0)(a\cdot\alpha+J_0)=a\cdot\alpha+J_0$ and $\alpha\cdot
a+J_0=\alpha\cdot((e+J_0)(a+J_0))=(\alpha\cdot
e+J_0)(a+J_0)=\alpha\cdot ea+J_0$, so
$e(a\cdot\alpha)-a\cdot\alpha \in J_0$ and $\alpha\cdot
ea-\alpha\cdot a \in J_0$, also $e(a\cdot\alpha)-(\alpha\cdot e)a
\in J\subseteq J_0$. Therefore $a\cdot\alpha-\alpha\cdot
a=a\cdot\alpha-e(a\cdot\alpha)+e(a\cdot\alpha)-\alpha\cdot
ea+\alpha\cdot ea-\alpha\cdot a \in J_0$. $\hfill\blacksquare $

We say the Banach algebra ${\mathfrak A}$ acts trivially on
$\mathcal A$ from left (right) if for each $\alpha\in \mathfrak A$
and $a\in \mathcal A$, $\alpha\cdot a=f(\alpha)a$
($a\cdot\alpha=f(\alpha)a$), where $f$ is a continuous linear
functional on ${\mathfrak A}$.

\vspace{.3cm}\paragraph{\large\bf Proposition 3.2.} {\it Let
$\mathcal A$ be module amenable as an ${\mathfrak A}$-module with
trivial left action, and let $J_0$ be a closed ideal of $\mathcal
A$ such that $J \subseteq J_0$. If $\mathcal A/J_0$ has an
identity, then $\mathcal A/J_0$ is amenable.}

\vspace{.3cm}\paragraph{\large\bf Proof.} Let $\mathcal X$ be a
unital $\mathcal A/J_0$-bimodule and $D: \mathcal
A/J_0\longrightarrow \mathcal X^* $ be a bounded derivation (see
\cite[Lemma 43.6]{bdu}). Then $\mathcal X$ is an $\mathcal
A$-bimodule with module actions given by
$$a\cdot x:=(a+J_0)\cdot x,\hspace{0.2cm}x\cdot a:=x\cdot(a+J_0) \hspace{0.3cm} (x \in \mathcal X ,a \in
\mathcal A),$$ and $\mathcal X$ is $\mathfrak A$-module with
trivial actions, that is $\alpha\cdot x=x\cdot\alpha=f(\alpha)x$,
for each $x\in \mathcal X$ and $\alpha \in \mathfrak A$ where $f$
is a continuous linear functional on ${\mathfrak A}$. Since
$f(\alpha)a-a\cdot\alpha\in J_0$, we have
$f(\alpha)a+J_0=a\cdot\alpha+J_0$, for each $\alpha\in \mathfrak
A$, and the actions of $\mathfrak A$ and $\mathcal A$ on
$\mathcal X$ are compatible. Therefore $\mathcal X$ is
commutative Banach ${\mathcal A}$-$\mathfrak A$-module. Consider
$\tilde{D}: \mathcal A \longrightarrow \mathcal X^*$ defined by
$\tilde{D}(a)=D(a+J_0)\hspace{0.1cm}(a \in \mathcal A)$. For $
a,b \in \mathcal A$ we have $\tilde{D}(a\pm b)=\tilde{D}(a)\pm
\tilde{D}(b)$ and $\tilde{D}(ab)=\tilde{D}(a)\cdot
b+a\cdot\tilde{D}(b).$ Also $\mathcal A/J_0$ is an $\mathfrak
A$-module, hence for $\alpha \in \mathfrak A$, we have

$$\tilde{D}(a\cdot\alpha)=D(a\cdot\alpha+J_0)=D(f(\alpha)a+J_0)=f(\alpha)D(e+J_0)=
D(e+J_0)\cdot\alpha.$$

On the other hand, since the left $\mathfrak A$-module actions on
$\mathcal A$ and $\mathcal X$ are trivial, $\tilde{D}(\alpha\cdot
a)=\tilde{D}(f(\alpha)a)=f(\alpha)D(a+J_0)=\alpha\cdot\tilde{D}(a)$.
Therefore there exists $x^* \in \mathcal X^*$ such that
$\tilde{D}(a)=a\cdot x^*-x^*\cdot a$ , hence
$D(a+J_0)=(a+J_0)\cdot x^*-x^*\cdot(a+J_0)$, and so $D$ is
inner.$\hfill\blacksquare $

\vspace{.3cm}\paragraph{\large\bf Proposition 3.3.} {\it If
${\mathfrak A}$ has a bounded approximate identity for $\mathcal
A$, then amenability of $\mathcal A/J$ implies module amenability
of $\mathcal A$.}

\vspace{.3cm}\paragraph{\large\bf Proof.} Let ${\mathcal X}$ be a
commutative Banach ${\mathcal A}$-${\mathfrak A}$-module. Since
$J\cdot{\mathcal X}={\mathcal X}\cdot J=0$, the following module
actions are well-defined
 $$(a+J)\cdot x:=a\cdot x, \hspace{0.2cm}x\cdot (a+J):=x\cdot a \hspace{0.3cm} (x \in \mathcal X ,a \in
\mathcal A),$$
 therefore ${\mathcal X}$ is a Banach $\mathcal A/J$-module.
 Suppose that $D: \mathcal A \longrightarrow {\mathcal X^*} $ is
a module derivation, and consider $\tilde{D}: \mathcal
A/J\longrightarrow \mathcal X^* $ defined by $\tilde{D}(a+J)=
 D(a)\, (a \in \mathcal A )$. We have
$$
D(\alpha\cdot ab-ab\cdot\alpha)=\alpha\cdot
D(ab)-D(ab)\cdot\alpha=0.
$$
 By the above observation, $\tilde{D}$ is
also well-defined. Since ${\mathfrak A}$ has a bonded approximate
identity for $\mathcal A$, it follows from the proof of
\cite[Proposition 2.1]{am1} that $\tilde{D}$ is
$\mathbb{C}$-linear, and so it is inner. Therefore $D$ is an
inner module derivation.$\hfill\blacksquare $

\vspace{.3cm} If $S$ is an inverse semigroup with an upward
directed set of idempotents $E$, then $E$ satisfies condition
$D_1$ of Duncan and Namioka, so $\ell ^{1}(E)$ has a bounded
approximate identity \cite{dna}. If $(\delta_{e_j})$ is a bounded
approximate identity of $\ell ^{1}(E)$, then $\delta_{e_j}*
\delta_s=\delta_{e_j}*
 \delta_{ss^*s}=\delta_{e_jss^*}*\delta_s \longrightarrow \delta_s$,
 and similarly for the right side multiplication. Therefore $\ell ^{1}(E)$ has a
 bounded approximate identity for
$\ell ^{1}(S)$. We use this fact to prove following result, which
is the main theorem of this section.

\vspace{.3cm}\paragraph{\large\bf Theorem 3.4.} {\it Let $S$ be an
inverse semigroup with an upward directed set of idempotents $E$.
 Then $\ell ^{1}(S)^{**}$ is module amenable (as an $\ell
^{1}(E)$-module with trivial left action) if and only if the
discrete group $S/\approx$ is finite.}

\vspace{.3cm}\paragraph{\large\bf Proof.} Let us move for a moment
to the general case and let $N$ be the closed ideal of ${\mathcal
A^{**}}$ generated by $(\alpha \cdot F)\square G-F \square
(G\cdot\alpha )$, for $F,G\in{\mathcal A}^{**}$ and $\alpha\in
{\mathfrak A}$. Then clearly $J\subseteq N$. Take two bounded nets
$(a_j),(b_k)\subset \mathcal A $ with $ \widehat a_j
\stackrel{J^{\perp}}{\longrightarrow}F$ and $ \widehat b_k
\stackrel{J^{\perp}}{\longrightarrow}G$, then
\begin{align*}
\langle{}(\alpha \cdot F)\square G-F \square (G\cdot\alpha
),f\rangle{} &=\langle{}\alpha .
(F\square G)-(F \square G)\cdot\alpha ,f\rangle{} \\
&=\langle{}F \square G,f\cdot\alpha - \alpha \cdot f\rangle{}\\
&=\lim_j\lim_k \langle{}\widehat a_j \widehat b_k ,f\cdot\alpha - \alpha \cdot f\rangle{}\\
&=\lim_j\lim_k \langle{}f,(\alpha \cdot a_j) b_k- a_j (b_k\cdot \alpha)\rangle{}\\
&=0.
\end{align*}
Therefore $ N \subseteq J^{\perp \perp}$. Hence we may consider
$J=N$ and $J_0=J^{\perp \perp}$ in Proposition 3.2 applied to
$\mathcal A^{**}$. Going back to the case where $\mathcal
A=\ell^1(S)$ and $\mathfrak A=\ell^1(E)$, since $S/\approx$ is a
discrete group, $\frac{\ell ^{1}(S)^{**}}{J^{\perp \perp}}\cong
\ell^1(S/\approx)^{**}$ has an identity. It follows from
Proposition 3.2 that $\ell^1(S/\approx)^{**}$ is amenable, hence
by \cite[Theorem 1.3]{glw}, $S/\approx$ is finite.

Conversely, if $S/\approx$ is finite, then $\ell ^{1}(S)/J$ is
amenable. Suppose that $\mathcal X$ is a Banach $\ell
^{1}(S)^{**}/N$-bimodule and $D:\ell
^{1}(S)^{**}/N\longrightarrow \mathcal X^*$ is a derivation. Then
$\mathcal X$ is $\ell ^{1}(S)/J$-bimodule with module actions
given by
$$(\delta_s+J)\cdot x:=(\delta_s+N)\cdot x,\hspace{0.3cm}x\cdot(\delta_s+J):=x\cdot(\delta_s+N),\hspace{0.3cm}
(x \in \mathcal X,s \in S).$$ Consider the linear map $\phi:\ell
^{1}(S)/J\longrightarrow \ell ^{1}(S)^{**}/N;\, \delta_s+J
\mapsto \delta_s+N$. It is clear that $D \circ \phi $ is a
derivation on the amenable Banach algebra $\ell ^{1}(S)/J$, and
so it is inner, hence $D$ is inner. Therefore $\ell
^{1}(S)^{**}/N$ is amenable. By Proposition 3.3, $ \ell
^1(S)^{**}$ is module amenable. $\hfill\blacksquare $

\vspace{.3cm} Consider the multiplication map $\omega_{\mathcal A}
\in {\mathfrak L}(\mathcal A \widehat \bigotimes \mathcal A,
\mathcal A )$ defined by $\omega_{\mathcal A} (a \otimes b)=ab$
and extended by linearity and continuity to a homomorphism of
Banach algebras. The module analog of this map is
$\widetilde{\omega}_{\mathcal A} \in {\mathfrak L}({\mathcal
A}\widehat \bigotimes _{\mathfrak A} {\mathcal A},\mathcal A/J)$
defined by $ \widetilde{\omega}_{\mathcal A} (a \otimes b
+I)=ab+J$, extended to a ${\mathfrak A}$-module homomorphism.
Then $ \widetilde{\omega}^*_{\mathcal A}$
 and $\widetilde{\omega}^{**}_{\mathcal A}$, the first and second adjoint of
 $ \widetilde{\omega}_{\mathcal A} $
  are ${\mathfrak A}$-module homomorphisms. It is proved in \cite[Lemma 1.7]{glw} that there is a continuous linear
mapping $\Omega :
 {\mathcal A}^{**}\widehat \bigotimes  {\mathcal
A}^{**} \longrightarrow ({\mathcal A}\widehat \bigotimes
 {\mathcal A})^{**}$ such that for $a,b,x \in
\mathcal A $ and $m \in {\mathcal A}^{**}\widehat \bigotimes
{\mathcal A}^{**}$ the following hold:

$(a)\,  \Omega (a \otimes b)= a \otimes b$,

$(b)\,  \Omega (m)\cdot x= \Omega(m\cdot x)$,

$(c)\,  x\cdot\Omega (m)= \Omega(x\cdot m)$,

$(d)\,  (\omega_{\mathcal A})^{**}(\Omega(m))= \omega_{\mathcal
A^{**}}(m)$.

We want to have a similar result for the module multiplication
map. To this end, let us first briefly go over the proof of the
above result. We have an isometric isomorphism between the space
of bilinear maps from ${\mathcal A}\times {\mathcal A}$ into
 $\mathbb{C}$ and $({\mathcal A}\widehat \bigotimes {\mathcal A})^*$
given by $T \mapsto \psi_T$ where $\psi_T(a \otimes b)= T(a,b).$
The map $ \Omega $ is then defined by
$$ \langle{}\Omega (F \otimes G),\psi_T\rangle{}= \lim_j\lim_k \langle{}\psi_T,a_j \otimes
b_k\rangle{},$$ where $(a_j),(b_k)$ are bounded nets in $
\mathcal A $
 such that  $ \widehat a_j \stackrel{w^*}{\longrightarrow}F$
and $ \widehat b_k \stackrel{w^*}{\longrightarrow}G$. Now let
${\mathcal A}^{**}\widehat \bigotimes _{\mathfrak A} {\mathcal
A}^{**}$ be the module projective  tensor product of ${\mathcal
A}^{**}$ and ${\mathcal A}^{**}$, that is ${\mathcal
A}^{**}\widehat \bigotimes _{\mathfrak A} {\mathcal
A}^{**}\cong\mathcal A^{**} \widehat \bigotimes \mathcal
A^{**}/M$, where $M$ is the closed ideal generated by elements of
the form $\alpha \cdot F \otimes G-F \otimes G\cdot\alpha$ for $
\alpha\in {\mathfrak A},$ and $ F,G\in{\mathcal A}^{**}$. Define
$ \Omega_{\mathfrak A}: \mathcal A^{**} \widehat \bigotimes
\mathcal A^{**}/M \longrightarrow ({\mathcal A}\widehat \bigotimes
 {\mathcal A})^{**}/I^{\perp \perp}$ via
 $$ \Omega_{\mathfrak A}(F \otimes G+M)= \Omega (F \otimes G)+ I^{\perp
 \perp}.$$
Then for $\varphi \in I^{\perp}$ and $\alpha \in {\mathfrak
  A}$ we have
 $$\langle{}\Omega(\alpha \cdot F \otimes G-F\otimes G\cdot\alpha),\varphi \rangle{}=
 \lim_j\lim_k \langle{}\varphi,\alpha\cdot a_j \otimes
b_k-a_j \otimes b_k\cdot\alpha\rangle{}= 0,$$ hence $
\Omega_{\mathfrak A}$ is well-define. Let $N$ be the closed ideal
defined in the proof of Theorem 3.4, then we know that $N
\subseteq J^{\perp \perp}$, so the map $\lambda :
A^{**}/N\longrightarrow A^{**}/J^{\perp \perp}; F+N
\longrightarrow F+ J^{\perp \perp}$ is a well defined continuous
homomorphism. For $a,b,x \in \mathcal A $ and $m \in {\mathcal
A}^{**}\widehat \bigotimes {\mathcal A}^{**}$ we have the
following equalities:

(1) $ \Omega_{\mathfrak A} (a \otimes b+M)= a \otimes b+I^{\perp
\perp}$,

(2) $ \Omega_{\mathfrak A}(m+M)\cdot x= \Omega_{\mathfrak
A}(m\cdot x+M)$,

(3) $ x\cdot\Omega_{\mathfrak A}(m+M)= \Omega_{\mathfrak
A}(x\cdot m+M)$,

(4) $ \tilde{\omega}_{\mathcal A}^{**}(\Omega_{\mathfrak
A}(m+M))=\lambda \circ \tilde{\omega}_{\mathcal A^{**}}(m+M)$.

Parts (1)-(3) are proved similar to $(a)$-$(c)$ in \cite[Lemma
1.7]{glw}. To see (4), given $\bar{a}=a+J \in \mathcal A/J, F$ and
$G$ in $\mathcal A^{**}$, take bounded nets $(a_j),(b_k)$ in $
\mathcal A $
 with $ \widehat a_j \stackrel{w^*}{\longrightarrow}F$
and $ \widehat b_k \stackrel{w^*}{\longrightarrow}G$, then
\begin{align*}
\langle{}\tilde{\omega}_{\mathcal A}^{**}(\Omega_{\mathfrak A}(F
\otimes G+M),\bar{a}\rangle{} &=\langle{}\Omega (F \otimes G)+
I^{\perp
 \perp},\tilde{\omega}_{\mathcal A}^{*}(\bar{a})\rangle{} \\
&=\lim_j\lim_k \langle{}\tilde{\omega}_{\mathcal A}^{*}(\bar{a}),a_j\otimes b_k+I)\rangle{}\\
&=\lim_j\lim_k \langle{}\bar{a},a_jb_k+J\rangle{}\\
&=\langle{}F \square G+J^{\perp \perp},\bar{a}\rangle{}\\
&=\langle{}\lambda(F \square G+N),\bar{a}\rangle{}\\
&=\langle{}\lambda \circ \tilde{\omega}_{\mathcal A^{**}}(F
\otimes G+M),\bar{a}\rangle{}.
\end{align*}

\vspace{.3cm}\paragraph{\large\bf Definition 3.5.} \cite{am1} {\it
A bounded net $\{\widetilde{\xi_j} \}$ in ${\mathcal A}\widehat
\bigotimes _{\mathfrak A} {\mathcal A}$ is called a {\it module
approximate diagonal} if $\widetilde{\omega}_{\mathcal A}(
\widetilde{\xi_j}) $ is a bounded approximate identity of  $
\mathcal A/J$ and
$$ \lim_j \| \xi_j \cdot a- a\cdot \xi_j \| =0 \hspace{0.3cm} (a \in \mathcal A
). $$ An element $ \widetilde{E} \in ({\mathcal A}\widehat
\bigotimes _{\mathfrak A} {\mathcal A})^{**}$ is called a {\it
module virtual diagonal} if
 $$ \widetilde{\omega}^{**}_{\mathcal A} (
\widetilde{E})\cdot a= \widetilde{a},\hspace{0.3cm}
\widetilde{E}\cdot a=a\cdot \widetilde{E} \hspace{0.3cm} (a \in
\mathcal A),
$$ where $ \widetilde{a}=a+ J^{\perp \perp }$.}

\vspace{.3cm}\paragraph{\large\bf Proposition 3.6.} {\it Let
$\mathcal A$ be a commutative Banach ${\mathfrak A}$-module such
that $\mathcal A^{**}$ is ${\mathfrak A}$-module amenable, then
so is $\mathcal A$.} \vspace{.3cm}\paragraph{\large\bf Proof.} We
have $J=\{0\}$, hence $N=J^{\perp\perp}=\{0\}$, in this case.
Since $\mathcal A^{**}$ is module amenable, by \cite[Proposition
2.2]{am1} it has a bounded approximate identity, so $\mathcal
A^{**}$ has a module approximate diagonal $\{\widetilde{m_i} \}$
in ${\mathcal A}^{**}\widehat \bigotimes _{\mathfrak A} {\mathcal
A}^{**}$ by \cite[Theorem 2.1]{am1}. For $F \in {\mathcal
A^{**}}$ and $ a \in \mathcal A $ we have $
\tilde{\omega}_{\mathcal A^{**}}(\widetilde{m_i}) {F}
\longrightarrow {F}$ and $\widetilde{m_i}\cdot
F-F\cdot\widetilde{m_i} \longrightarrow 0$. By an arguments
similar to what is used to prove part (4) above, we get
$\tilde{\omega}_{\mathcal A}^{**}(\Omega_{\mathfrak
A}(\widetilde{m_i})){F} \longrightarrow {F}$ and
$\Omega_{\mathfrak A}(\widetilde{m_i})\cdot
a-a\cdot\Omega_{\mathfrak A}(\widetilde{m_i}) \longrightarrow 0
$. Since $\{\Omega_{\mathfrak A}(\widetilde{m_i}) \}$ is bounded,
it has a cluster point $ \widetilde{E}$ in $({\mathcal A}\widehat
\bigotimes
 {\mathcal A})^{**}/I^{\perp \perp}\cong ({\mathcal A}\widehat
\bigotimes _{\mathfrak A} {\mathcal A})^{**}$, such that
 $\tilde{\omega}_{\mathcal A}^{**}(\widetilde{E})
 \cdot a=\tilde{a}$ and $\widetilde{E}\cdot a=a\cdot \widetilde{E}$. Therefore
 $\widetilde{E}$ is a module virtual diagonal for $\mathcal A$, and by \cite[Theorem 2.1]{am1} $\mathcal A$ is module amenable.$\hfill\blacksquare $

\vspace{.3cm} We finish by two examples. First we give an example
for which $\mathcal A$ is module amenable but $\mathcal A^{**}$
is not. In the second example $\mathcal A^{**}$ is module
amenable but it is not amenable. Note that in these examples
$\mathcal A$ is not a commutative $\mathfrak A$-module, as the
left action is taken to be trivial. Let $\mathcal C$ be the
bicyclic inverse semigroup generated by $p$ and $q$, that is
$$\mathcal C=\{p^mq^n : m,n\geq 0 \},\hspace{0.2cm}(p^mq^n)^*=p^nq^m. $$ The
multiplication operation is defined by
$$(p^mq^n)(p^{m'}q^{n'})=p^{m-n+max \{n,m'\}}q^{m'-n'+max \{n,m'\}}.$$
The set of idempotents of $\mathcal C$ is $E_{\mathcal
C}=\{p^nq^n : n=0,1,...\}$ which is totally ordered with the
following order
$$p^nq^n \leq p^mq^m \Longleftrightarrow m \leq n.$$
 $\mathcal C$ is isomorphic
to the semigroup $$\langle{}\{e,t,s \mid se=es=s, te=et=t, st=e,
ts\neq e \}\rangle{}.$$ Consider the equivalence relation
$\approx$ on $\mathcal C$ defined before Theorem 2.4. Clearly
$es\approx se\approx s $. Also since $ts$ and $st$ are
idempotents, $ts\approx st\approx e $. Therefore
$$\mathcal C/\approx=\langle{}\{s,t,e |st \approx ts \approx e, se\approx
es\approx s\}\rangle{}$$ which is the cyclic group generated by
$s$. The element $s$ is not an idempotent, so $s$ is not
equivalent to $e$. If for some $k \in \mathbb{N}, s^k\approx
s^{k-1}$, then $ \delta_{s^k}-\delta_{s^{k-1}} \in J$, so
$\delta_s-\delta_e=\delta_{s^kt^{k-1}}-\delta_{s^{k-1}t^{k-1}} \in
J$, therefore $s\approx e$, which is a contradiction. Thus
$\mathcal C/\approx$ is isomorphic to $\mathbb{Z}$, and hence
${\ell ^{1}}(\mathcal C/\approx)$ is amenable. It follows from
Theorem 3.3 that ${\ell ^{1}}(\mathcal C)$ is ${\ell
^{1}}(E_{\mathcal C})$-module amenable. On the other hand, since
$\mathcal C/\approx$ is infinite, by Theorem 3.4, ${\ell
^{1}}(\mathcal C)^{**}$ is not ${\ell ^{1}}(E_{\mathcal
C})$-module amenable.

Next let $(\mathbb{N}, \vee)$ be the semigroup of positive
integers with maximum operation, that is $m\vee n=max(m,n)$,
then  each element of $\mathbb{N}$ is an idempotent, hence
$\mathbb{N}/\approx$ is the trivial group with one element.
Therefore ${\ell ^{1}}( \mathbb{N})^{**}$ is ${\ell ^{1}}(
E_{\mathbb{N}})$-module amenable. Since $\mathbb{N}$ is infinite
weakly cancellative semigroup, ${\ell ^{1}}( \mathbb{N})^{**}$ is
not amenable \cite[Theorem 1.3]{glw}.

\end{document}